\documentclass{amsart}
\usepackage{amsmath,amssymb,amsbsy,amsfonts,amsthm,latexsym,
amsopn,amstext,amsxtra,euscript,amscd,graphicx}
\newtheorem{teor}{Theorem}

\theoremstyle{definition}

\begin{document}

\title{Further generalizations of the parallelogram law}

\author{Antonio M. Oller-Marc\'{e}n}
\address{Centro Universitario de la Defensa de Zaragoza\\ Ctra. Huesca s/n, 50090 Zaragoza, Spain} 
\email{oller@unizar.es}

\maketitle

\begin{abstract}
In recent work \cite{FON}, a generalization of the parallelogram law in any dimension $N\geq 2$ was given by considering the ratio of the quadratic mean of the measures of the $N-1$-dimensional diagonals to the quadratic mean of the measures of the faces of a parallelotope. In this paper, we provide a further generalization considering not only $(N-1)$-dimensional diagonals and faces, but the $k$-dimensional ones for every $1\leq k\leq N-1$.
\end{abstract}

\section{Introduction}
If we consider the usual Euclidean space $(\mathbb{R}^{n},\|\cdot\|)$, the well-known identity
\begin{equation}\label{ol}
\|a+b\|^2+\|a-b\|^2=2(\|a\|^2+\|b\|^2)
\end{equation}
is called the \textit{parallelogram law}.

This identity can be extended to higher dimensions in several ways. For example, it is straightforward to see that
\begin{equation}\label{genl}
\|a+b+c\|^2+\|a+b-c\|^2+\|a-b+c\|^2+\|a-b-c\|^2=4(\|a\|^2+\|b\|^2+\|c\|^2)
\end{equation}
with subsequent analogue identities arising inductively. There are, in fact, many works devoted to provide generalizations of (\ref{ol}) in many different contexts \cite{ome,nash,pen}

Note that if we rewrite (\ref{ol}) as
\begin{equation}
\frac{\|a+b\|^2+\|a-b\|^2}{2}=2\frac{(\|a\|^2+\|b\|^2+\|a\|^2+\|b\|^2)}{4}
\end{equation}
it just means that in any parallelogram, the ratio of the quadratic mean of the lengths of its diagonals to the quadratic mean of the lengths of its sides equals $\sqrt{2}$. With this interpretation in mind, Alessandro Fonda \cite{FON} has recently proved the following interesting generalization.

\begin{teor}\label{FON}
Given linearly independent vectors $a_1,\dots, a_N\in\mathbb{R}^n$, it holds that
\begin{align*}
&\sum_{i<j}\left(\left\|(a_i+a_j)\wedge \bigwedge_{k\neq i,j} a_k\right\|^2+\left\|(a_i-a_j)\wedge \bigwedge_{k\neq i,j} a_k\right\|^2\right)=\\&=(N-1)\sum_{k=1}^{N}2\left\|a_1\wedge\cdots\wedge\widehat{a_k}\wedge\cdots\wedge a_N\right\|^2.
\end{align*}
In other words, for any $N$-dimensional parallelotope, the ratio of the quadratic mean of the $(N-1)$-dimensional measures of its diagonals to the quadratic mean of the $(N-1)$-dimensional measures of its faces is equal to $\sqrt{2}$. 
\end{teor}

In this work we extend this result considering the faces of dimension $k$ for every $1\leq k\leq N-1$ and a suitable definition of $k$-dimensional diagonal of a parallelotope. Then, Theorem \ref{FON} will just be a particular case of our result for $k=N-1$. Indeed, our result can be stated as follows.

\begin{teor}\label{mio}
Let us consider an $N$-dimensional parallelotope and let $1\leq k\leq N-1$. The ratio of the quadratic mean of the $k$-dimensional measures of its $k$-dimensional diagonals to the quadratic mean of the $k$-dimensional measures of its $k$-dimensional faces is equal to $\sqrt{N-k+1}$. 
\end{teor}

In fact, our generalization goes in the line of the work \cite{nash} but considering the definition of diagonal face given in \cite{FON}. 

\section{Notation and preliminaries}

In this section we are going to introduce some notation and to present some basic facts that will be useful in the sequel. Let us consider a parallelotope $\mathcal{P}$ generated by a family of linearly independent vectors $\mathcal{B}=\{a_1,a_2,\dots,a_N\}\subseteq\mathbb{R}^n$. This means that
$$\mathcal{P}=\left\{\sum_{i=1}^N\alpha_i a_i:\alpha_i\in[0,1]\right\}.$$

Let us fix $1\leq k\leq N-1$. Then, given $k$ different vectors $\mathcal{S}=\{a_{i_1},\dots,a_{i_k}\}\subseteq\mathcal{B}$, we can consider the face generated by them:
$$\mathcal{F}(\mathcal{S})=\left\{\sum_{v\in\mathcal{S}}\alpha_v v:\alpha_v\in[0,1]\right\}.$$
This face can now be translated by one or more of the remaining vectors thus obtaining a face 
$$\mathcal{F}^{I}(\mathcal{S})=\left\{\sum_{v\in\mathcal{S}}\alpha_v a_{v}+\sum_{w\in\mathcal{B}\setminus\mathcal{S}}\alpha_w w\in\mathcal{P}:\alpha_w\in\{0,1\}\right\},$$
where $I=(\alpha_v)_{v\not\in\mathcal{S}}\in\{0,1\}^{N-k}$. Since each choice of a set $\mathcal{S}\subseteq\mathcal{B}$ and a vector $I\in\{0,1\}^{N-k}$ leads to a different face and every face can be obtained in this way, it follows the well-known result that $\mathcal{P}$ has exactly $2^{N-k}\binom{N}{k}$ $k$-dimensional faces. Moreover, it is clear that all the $2^{N-k}$ different faces $\mathcal{F}^{I}(\mathcal{S})$ are congruent to the set generated by $\mathcal{S}$, $\mathcal{F}(\mathcal{S})$.

Now, we focus on the $k$-dimensional diagonals which will be defined following the ideas in \cite{FON}. Let us consider $N-k+1$ different vectors $\mathcal{T}=\{a_{i_1},\dots,a_{i_{N-k+1}}\}\subseteq\mathcal{B}$ and let $\mathcal{T}=\mathcal{T}_1\cup\mathcal{T}_2$ be any partition. Then, the following set
$$\mathcal{D}(\mathcal{T}_1,\mathcal{T}_2)=\left\{\alpha\sum_{v\in\mathcal{T}_1}v+(1-\alpha)\sum_{v\in\mathcal{T}_2}v+\sum_{w\in\mathcal{B}\setminus\mathcal{T}}\alpha_w w:\alpha,\alpha_{w}\in[0,1]\right\}.$$
is called the $k$-dimensional diagonal associated to $(\mathcal{T},\mathcal{T}_1,\mathcal{T}_2)$. Clearly each choice of a set $\mathcal{T}\subseteq\mathcal{B}$ and a partition of $\mathcal{T}$ leads to a different diagonal. Thus, it readily follows that $\mathcal{P}$ has exactly $2^{N-k}\binom{N}{N-k+1}$ different $k$-dimensional diagonals. Moreover, if we define the vector
$$V(\mathcal{T}_1,\mathcal{T}_2)=\sum_{v\in\mathcal{T}_1}v-\sum_{v\in\mathcal{T}_2}v,$$
we have that 
$$\mathcal{D}(\mathcal{T}_1,\mathcal{T}_2)=\left\{\alpha V(\mathcal{T}_1,\mathcal{T}_2)+\sum_{v\in\mathcal{T}_2}v+\sum_{w\in\mathcal{B}\setminus\mathcal{T}}\alpha_w w:\alpha,\alpha_w\in[0,1]\right\}$$
and, consequently, it is clear that the diagonal $\mathcal{D}(\mathcal{T}_1,\mathcal{T}_2)$ is just a translation of the set generated by $\{V(\mathcal{T}_1,\mathcal{T}_2),w:w\in\mathcal{B}\setminus\mathcal{T}\}$ and, hence, it is congruent to it.

\section{Proof of Theorem \ref{mio}}

After introducing the notation and the main objects involved in thie work, we are now in the condition to proof the main result of the paper. 

Let $\mathcal{P}$ be a parallelotope generated by $\mathcal{B}=\{a_1,a_2,\dots,a_N\}\subseteq\mathbb{R}^n$. We first compute the quadratic mean of the $k$-dimensional measures of its $k$-dimensional faces. To do so, we first note that, for every $\mathcal{S}=\{a_{i_1},\dots,a_{i_k}\}\subseteq\mathcal{B}$, the $k$-dimensional measure of the face $\mathcal{F}(\mathcal{S})$ is $\|a_{i_1}\wedge\dots\wedge a_{i_k}\|$. In the previous section we have seen that $\mathcal{P}$ has exactly $2^{N-k}\binom{N}{k}$ $k$-dimensional faces and, moreover, that there are exactly $2^{N-k}$ copies of each face $\mathcal{F}(\mathcal{S})$. Consequently, the quadratic mean of the $k$-dimensional measures of the $k$-dimensional faces of $\mathcal{P}$ is:
\begin{equation}\label{qf}\frac{\displaystyle 2^{N-k}\sum\|a_{i_1}\wedge\dots\wedge a_{i_k}\|^2}{\displaystyle 2^{N-k}\binom{N}{k}}.\end{equation}

Now we have to compute the quadratic mean of the $k$-dimensional measures of the $k$-dimensional diagonals of $\mathcal{P}$. First of all, recall that $\mathcal{P}$ has exactly $2^{N-k}\binom{N}{N-k+1}$ different $k$-dimensional diagonals. Each of them is the translation of the set generated by $\{V(\mathcal{T}_1,\mathcal{T}_2),w:w\in\mathcal{B}\setminus\mathcal{T}\}$ for exactly one choice of $(\mathcal{T},\mathcal{T}_1,\mathcal{T}_2)$. The $k$ dimensional measure of this latter set is $\displaystyle \left\|V(\mathcal{T}_1,\mathcal{T}_2)\wedge\bigwedge_{w\in\mathcal{B}\setminus\mathcal{T}} w\right\|$. Consequently, the quadratic mean of the $k$-dimensional measures of the $k$-dimensional diagonals of $\mathcal{P}$ is:
\begin{equation}\label{qd}\frac{\displaystyle \sum_{\mathcal{T},\mathcal{T}_1,\mathcal{T}_2}\left\|V(\mathcal{T}_1,\mathcal{T}_2)\wedge\bigwedge_{w\in\mathcal{B}\setminus\mathcal{T}} w\right\|^2}{\displaystyle 2^{N-k}\binom{N}{N-k+1}}.\end{equation}
Now, using the bilinearity of the scalar product and taking into account the definition of $V(\mathcal{T}_1,\mathcal{T}_2)$, it can be easily seen that when we vary $(\mathcal{T},\mathcal{T}_1,\mathcal{T}_2)$, we get the term $\|a_{i_1}\wedge\dots\wedge a_{i_k}\|^2$ exactly $2^{N-K}k$ times for every possible choice of $\{a_{i_1},\dots,a_{i_k}\}\subseteq\mathcal{B}$. This implies that the quadratic mean of the $k$-dimensional measures of the $k$-dimensional diagonals of $\mathcal{P}$ (\ref{qd}) can in fact be written as:
\begin{equation}\label{qdf}\frac{\displaystyle 2^{N-k}k\sum\|a_{i_1}\wedge\dots\wedge a_{i_k}\|^2}{\displaystyle 2^{N-k}\binom{N}{N-k+1}}.\end{equation}

Finally, in order to obtain Theorem \ref{mio} it is enough to divide (\ref{qdf}) by (\ref{qf}):
$$\frac{(\ref{qdf})}{(\ref{qf})}=\frac{k\binom{N}{k}}{\binom{N}{N-k+1}}=N-k+1.$$

\end{document}